\input amstex
\documentstyle{amsppt}
\nologo
\pageno=1
\loadbold
\leftheadtext{{\smc Nandor Simanyi}}
\rightheadtext{{\smc Proving The Ergodic Hypothesis...}}

\hsize=5 true in
\vsize=8.4 true in
\hoffset=.75 in
\TagsOnRight
%\NoBlackBoxes

\define\flow{\left(\bold{M},\{S^t\}_{t\in\Bbb R},\mu\right)}
\define\symb{\Sigma=\left(\sigma_1,\sigma_2,\dots,\sigma_n\right)}

\define\endrem{}

\footnote""{AMS 1991 Subject Classification: Primary 37D50, 34D05.}
\footnote""{Key words and phrases: Non-uniformly hyperbolic dynamical systems,
singularities, cylindric billiards, ergodicity, mixing, local ergodicity.}

\footnote""{This research was partially supported by the National Science 
Foundation, grant DMS-0098773.}

%%\email{ h-kubo@@math.hokudai.ac.jp}
\bigskip
\document

\vglue 1\baselineskip

\centerline{\bf  Proving The Ergodic Hypothesis for Billiards With Disjoint 
Cylindric Scatterers}

%\centerline{\bf  }

\bigskip

\medskip

\centerline{\smc Nandor Simanyi}
\medskip

{\eightpoint
\centerline{ University of Alabama at Birmingham}
\centerline{ Department of Mathematics}
\centerline{ Campbell Hall, Birmingham, AL 35294 U.S.A.}
}
\medskip

{\eightpoint{\narrower\smallskip\noindent {\bf Abstract.}
We study the ergodic properties of mathematical billiards describing the 
uniform motion of a point in a flat torus from which finitely many, pairwise
disjoint, tubular neighborhoods of translated subtori (the so called cylindric
scatterers) have been removed. We prove that every such system is ergodic
(actually, a Bernoulli flow), unless a simple geometric obstacle for the
ergodicity is present.
\smallskip}}
\medskip

\noindent {\bf 1. Introduction}
Non-uniformly hyperbolic systems (possibly, with singularities) play a pivotal
role in the ergodic theory of dynamical systems. Their systematic study
started several decades ago, and it is not our goal here to provide the reader
with a comprehensive review of the history of these investigations but,
instead, we opt for presenting in nutshell a cross section of a few selected
results.

In 1939 G. A. Hedlund and E. Hopf [He(1939)], [Ho(1939)], proved the
hyperbolic ergodicity of geodesic flows on closed, compact surfaces with
constant negative curvature by inventing the famous method of "Hopf chains"
constituted by local stable and unstable invariant manifolds.

In 1963 Ya. G. Sinai [Sin(1963)] formulated a modern version of Boltzmann's
ergodic hypothesis, what we call now the "Boltzmann-Sinai ergodic hypothesis":
the billiard system of $N$ ($\ge2$) hard spheres of unit mass moving in the
flat torus $\Bbb T^\nu=\Bbb R^\nu/\Bbb Z^\nu$ ($\nu\ge2$) is ergodic after
we make the standard reductions by fixing the values of the trivial invariant
quantities. It took seven years until he proved this conjecture for the
case $N=2$, $\nu=2$ in [Sin(1970)]. Another 17 years later N. I. Chernov
and Ya. G. Sinai [S-Ch(1987)] proved the hypothesis for the case $N=2$,
$\nu\ge2$ by also proving a powerful and very useful theorem on local
ergodicity.

In the meantime, in 1977, Ya. Pesin [P(1977)] laid down the foundations
of his theory on the ergodic properties of smooth, hyperbolic dynamical
systems. Later on this theory (nowadays called Pesin theory) was
significantly extended by A. Katok and J-M. Strelcyn [K-S(1986)]
to hyperbolic systems with singularities. That theory is already
applicable for billiard systems, too.

Until the end of the seventies the phenomenon of hyperbolicity (exponential
instability of the trajectories) was almost exclusively attributed to some
direct geometric scattering effect, like negative curvature of space, or
strict convexity of the scatterers. This explains the profound shock that
was caused by the discovery of L. A. Bunimovich [B(1979)]: certain focusing
billiard tables (like the celebrated stadium) can also produce complete
hyperbolicity and, in that way, ergodicity. It was partly this result that
led to Wojtkowski's theory of invariant cone fields, [W(1985)], 
[W(1986)].

The big difference between the system of two spheres in $\Bbb T^\nu$
($\nu\ge2$, [S-Ch(1987)]) and the system of $N$ ($\ge3$) spheres in
$\Bbb T^\nu$ is that the latter one is merely a so called semi-dispersive
billiard system (the scatterers are convex but not strictly convex
sets, namely cylinders), while the former one is strictly dispersive
(the scatterers are strictly convex sets). This fact makes the proof
of ergodicity (mixing properties) much more complicated. In our series
of papers jointly written with A. Kr\'amli and D. Sz\'asz [K-S-Sz(1990)],
[K-S-Sz(1991)], and [K-S-Sz(1992)], we managed to prove the (hyperbolic)
ergodicity of three and four billiard spheres in the toroidal container
$\Bbb T^\nu$. By inventing new topological methods and the Connecting Path
Formula (CPF), in my two-part paper [Sim(1992)] I proved the (hyperbolic)
ergodicity of $N$ hard spheres in $\Bbb T^\nu$, provided that $N\le\nu$.

The common feature of hard sphere systems is --- as D. Sz\'asz pointed this
out first in [Sz(1993)] and [Sz(1994)] --- that all of them belong to the
family of so called cylindric billiards, the definition of which can be
found later in this paragraph. However, the first appearance of a special,
3-D cylindric billiard system took place in [K-S-Sz(1989)], where we
proved the ergodicity of a 3-D billiard flow with two orthogonal
cylindric scatterers. Later D. Sz\'asz [Sz(1994)] presented a complete
picture (as far as ergodicity is concerned) of cylindric billiards with
cylinders whose generator subspaces are spanned by mutually orthogonal
coordinate axes. The task of proving ergodicity for the first non-trivial,
non-orthogonal cylindric billiard system was taken up in [S-Sz(1994)].

Finally, in our joint venture with D. Sz\'asz [S-Sz(1999)] we managed to
prove the complete hyperbolicity of {\it typical} hard sphere systems.

\subheading{\bf Cylindric billiards} Consider the $d$-dimensional
($d\ge2$) flat torus $\Bbb T^d=\Bbb R^d/\Cal L$ supplied with the
usual Riemannian inner product $\langle\, .\, ,\, .\, \rangle$ inherited
from the standard inner product of the universal covering space $\Bbb R^d$.
Here $\Cal L\subset\Bbb R^d$ is supposed to be a lattice, i. e. a discrete
subgroup of the additive group $\Bbb R^d$ with $\text{rank}(\Cal L)=d$.
The reason why we want to allow general lattices, other than just the
integer lattice $\Bbb Z^d$, is that otherwise the hard sphere systems would
not be covered. The geometry of the structure lattice $\Cal L$ in the
case of a hard sphere system is significantly different from the geometry
of the standard lattice $\Bbb Z^d$ in the standard Euclidean space
$\Bbb R^d$, see subsection 2.4 of [Sim(2002)].

The configuration space of a cylindric billiard is
$\bold Q=\Bbb T^d\setminus\left(C_1\cup\dots\cup C_k\right)$, where the
cylindric scatterers $C_i$ ($i=1,\dots,k$) are defined as follows:

Let $A_i\subset\Bbb R^d$ be a so called lattice subspace of $\Bbb R^d$,
which means that $\text{rank}(A_i\cap\Cal L)=\text{dim}A_i$. In this case
the factor $A_i/(A_i\cap\Cal L)$ is a subtorus in $\Bbb T^d=\Bbb R^d/\Cal L$
which will be taken as the generator of the cylinder 
$C_i\subset\Bbb T^d$, $i=1,\dots,k$. Denote by $L_i=A_i^\perp$ the
orthocomplement of $A_i$ in $\Bbb R^d$. Throughout this paper we will
always assume that $\text{dim}L_i\ge2$. Let, furthermore, the numbers
$r_i>0$ (the radii of the spherical cylinders $C_i$) and some translation
vectors $t_i\in\Bbb T^d=\Bbb R^d/\Cal L$ be given. The translation
vectors $t_i$ play a crucial role in positioning the cylinders $C_i$
in the ambient torus $\Bbb T^d$. Set

$$
C_i=\left\{x\in\Bbb T^d:\; \text{dist}\left(x-t_i,A_i/(A_i\cap\Cal L)
\right)<r_i \right\}.
$$
In order to avoid further unnecessary complications, we always assume that
the interior of the configuration space 
$\bold Q=\Bbb T^d\setminus\left(C_1\cup\dots\cup C_k\right)$ is connected.
The phase space $\bold M$ of our cylindric billiard flow will be the
unit tangent bundle of $\bold Q$ (modulo some natural identification
at its boundary), i. e. $\bold M=\bold Q\times\Bbb S^{d-1}$. (Here 
$\Bbb S^{d-1}$ denotes the unit sphere of $\Bbb R^d$.)

The dynamical system $\flow$, where $S^t$ ($t\in\Bbb R$) is the dynamics 
defined by uniform motion inside the domain $\bold Q$ and specular
reflections at its boundary (at the scatterers), and $\mu$ is the
Liouville measure, is called the cylindric billiard flow we want to
investigate. (As to the notions and notations in connection with 
semi-dispersive billiards, the reader is kindly recommended to consult
the work [K-S-Sz(1990)].)

\medskip

\subheading{\bf Transitive cylindric billiards}

The main conjecture concerning the (hyperbolic) ergodicity of cylindric
billiards is the "Erd\H otarcsa conjecture" (named after the picturesque
village in rural Hungary where it was originally formulated) that appeared
as Conjecture 1 in Section 3 of [S-Sz(2000)].

\medskip

\subheading{\bf Transitivity} Let $L_1,\dots,L_k\subset\Bbb R^d$
be subspaces, $\text{dim}L_i\ge2$, $A_i=L_i^\perp$, $i=1,\dots,k$. Set

$$
\Cal G_i=\left\{U\in\text{SO}(d):\, U|A_i=\text{Id}_{A_i}\right\},
$$
and let 
$\Cal G=\left\langle\Cal G_1,\dots,\Cal G_k\right\rangle\subset\text{SO}(d)$
be the algebraic generate of the compact, connected Lie subgroups
$\Cal G_i$ in $\text{SO}(d)$. The following notions appeared in Section 3 of
[S-Sz(2000)].

\subheading{\bf Definition} We say that the system of base spaces
$\{L_1,\dots,L_k\}$ (or, equivalently, the cylindric billiard system defined
by them) is {\it transitive} if and only if the group $\Cal G$ acts 
transitively on the unit sphere $\Bbb S^{d-1}$ of $\Bbb R^d$.

\subheading{Definition} We say that the system of subspaces
$\{L_1,\dots,L_k\}$ has the Orthogonal Non-splitting Property (ONSP) if there
is no non-trivial orthogonal splitting $\Bbb R^d=B_1\oplus B_2$ of
$\Bbb R^d$ with the property that for every index $i$ ($1\le i\le k$)
$L_i\subset B_1$ or $L_i\subset B_2$.

The next result can be found in Section 3 of [S-Sz(2000)] (see 3.1--3.6
thereof):

\subheading{\bf Proposition} For the system of subspaces
$\{L_1,\dots,L_k\}$ the following three properties are equivalent:

\medskip

(1) $\{L_1,\dots,L_k\}$ is transitive;

\medskip

(2) the system of subspaces $\{L_1,\dots,L_k\}$ has the ONSP;

\medskip

(3) the natural representation of $\Cal G$ in $\Bbb R^d$ is irreducible.

\medskip

\subheading{\bf The Erd\H otarcsa conjecture} A cylindric billiard flow is
ergodic if and only if it is transitive. In that case the cylindric billiard 
system is actually a completely hyperbolic Bernoulli flow, see [C-H(1996)] 
and [O-W(1998)].

\medskip

In order to avoid unnecessary complications, throughout the paper we always 
assume that

$$
\text{int}\bold Q\text{ is connected, and}
\tag 1.1
$$

$$
\aligned
\text{the }d\text{-dim spatial angle }\alpha(q)\text{ subtended by }\bold Q \\
\text{at any of its boundary points }q\in\partial\bold Q\text{ is positive.}
\endaligned
\tag 1.2
$$

\medskip

The Erd\H otarcsa Conjecture has not been proved so far in full generality.
Certain partial results, however, exist. Without pursuing the goal of
achieving completeness, here we cite just two of such results:

\medskip

\subheading{\bf Theorem of [Sim(2001)]} Almost every hard disk system (i. e.
hard sphere system in a $2$-D torus) is hyperbolic and ergodic. (Here 
``almost every'' is meant with respect to the outer geometric parameters
$(r;\, m_1,\dots,m_N)$, where $r>0$ is the common radius of the disks, while
$m_1,\dots,m_N$ are the masses.) 

\medskip

\subheading{\bf Theorem of [Sim(2002)]} Every hard sphere system is completely
hyperbolic, i. e. all of its relevant Lyapunov exponents are nonzero almost
everywhere.

\medskip

In this paper we are mainly interested in understanding the ergodic properties
of cylindric billiard flows $\flow$ in which the closures $\bar C_i$ of the
scattering cylinders $C_i$ are pairwise disjoint, i. e.

$$
\bar C_i\cap\bar C_j=\emptyset\text{ for }1\le i<j\le k.
\tag 1.3
$$
Elementary linear algebra shows that for such a (disjoint) cylindric billiard
system it is true that $\text{span}\left\{A_i,A_j\right\}\ne\Bbb R^d$
for $1\le i,\, j\le k$ or, equivalently, 

$$
L_i\cap L_j\ne\{0\}\text{ for }1\le i,\, j\le k.
\tag 1.4
$$
From now on we drop the disjointness condition (1.3) by only keeping the
somewhat relaxed condition (1.4) above. The first, very simple, question that
arises here is to characterize the transitivity of the $\Cal G$-action on
$\Bbb S^{d-1}$ under the condition (1.4) above. The following proposition
immediately follows from (1.4) and the characterization (2) of the 
transitivity above:

\medskip

\subheading{\bf Proposition 1.5} A cylindric billiard system with the 
additional property (1.4) is transitive (that is, the $\Cal G$-action on
the velocity sphere $\Bbb S^{d-1}$ is transitive) if and only if 
$\text{span}\left\{L_1,\dots,L_k\right\}=\Bbb R^d$ or, equivalently,
$\bigcap_{i=1}^k A_i=\{0\}$. \qed

\medskip

Now we are able to put forward the result of this paper:

\medskip

\subheading{\bf Theorem} Assume that the cylindric billiard flow $\flow$
enjoys the geometric properties (1.1), (1.2), and (1.4) above. Then the
transitivity condition $\bigcap_{i=1}^k A_i=\{0\}$ implies that the flow
$\flow$ is completely hyperbolic and ergodic.

\medskip

\subheading{\bf Remark} In the case of a hard sphere system with masses
$m_1,\,m_2,\dots,\,m_N$ the base space $L_{ij}$ of the cylinder $C_{ij}$
(describing the collision between the spheres labelled by $i$ and $j$)
is obviously the set

$$
\left\{(\delta q_1,\dots,\delta q_N)\in\Cal T\bold Q\big|\; \delta q_k=0
\text{ for } k\not\in\{i,\,j\}\right\}.
$$
(See also (4.4) in [S-Sz(2000)].) Therefore, in such systems the intersection
of any two base spaces is zero. This shows that our present result is
complementary to any possible past and future result about hard sphere systems.

\medskip

\subheading{\bf Organizing the Paper} Section 2 contains the indispensable
technical preparations, definitions, and notations. \S3 is devoted to
proving that if a non-singular orbit segment $S^{[a,b]}x$ of the flow
$\flow$ has a combinatorially rich symbolic collision sequence (in a well
defined sense) then $S^{[a,b]}x$ is sufficient (geometrically hyperbolic)
modulo a codimension-2 algebraic subset of the phase space. Finally, the
closing Section 4 contains the inductive proof of

\medskip

(H1) the so called ``Chernov--Sinai Ansatz'' for the flow $\flow$, i. e. that
--- informally speaking --- for almost every singular phase point
$x\in\bold M$ the forward semi-trajectory after the singularity is sufficient

\medskip

\noindent
and

\medskip

(H2) outside of a slim subset $S\subset\bold M$ (for the notion of slimness,
please see \S2 below) it is true that

(i) $S^{(-\infty,\infty)}x$ has at most one singularity;

(ii) $S^{(-\infty,\infty)}x$ is sufficient.

\medskip

Section 4 concludes with putting together all the above results and applying 
the Theorem on Local Ergodicity for semi-dispersive billiards
[S-Ch(1987)] to complete the proof of our Theorem.

\medskip

\subheading{\bf Remark} In order to simplify the notations, throughout the
paper we will assume that the fundamental lattice $\Cal L\subset\Bbb R^d$
of the factorization $\Bbb T^d=\Bbb R^d/\Cal L$ is the standard integer
lattice $\Bbb Z^d$. This is not a significant restriction of generality, for
the entire proof of the theorem carries over to the general case easily
by an almost word-by-word translation.

\medskip

\noindent {\bf 2. Prerequisites}
\subheading{Trajectory Branches}
We are going to briefly describe the discontinuity of the flow
$\{S^t\}$ caused by a multiple collision at time $t_0$.
Assume first that the pre--collision velocities of the particles are given.
What can we say about the possible post--collision velocities? Let us perturb
the pre--collision phase point (at time $t_0-0$) infinitesimally, so that the
collisions at $\sim t_0$ occur at infinitesimally different moments. By
applying the collision laws to the arising finite sequence of collisions, we
see that the post-collision velocities are fully determined by the time--
ordering of the considered collisions. Therefore, the collection of all
possible time-orderings of these collisions gives rise to a finite family of
continuations of the trajectory beyond $t_0$. They are called the
{\bf trajectory branches}. It is quite clear that similar statements can be
said regarding the evolution of a trajectory through a multiple collision
{\bf in reverse time}. Furthermore, it is also obvious that for any given
phase point $x_0\in\bold M$ there are two, $\omega$-high trees
$\Cal T_+$ and $\Cal T_-$ such that $\Cal T_+$ ($\Cal T_-$) describes all the
possible continuations of the positive (negative) trajectory
$S^{[0,\infty)}x_0$ ($S^{(-\infty,0]}x_0$). (For the definitions of trees and
for some of their applications to billiards, cf. the beginning of \S 5
in [K-S-Sz(1992)].) It is also clear that all possible continuations
(branches) of the whole trajectory $S^{(-\infty,\infty)}x_0$ can be uniquely
described by all possible pairs $(B_-,B_+)$ of $\omega$-high branches of
the trees $\Cal T_-$ and $\Cal T_+$ ($B_-\subset\Cal T_-, B_+\subset
\Cal T_+$).

Finally, we note that the trajectory of the phase point $x_0$ has exactly two
branches, provided that $S^tx_0$ hits a singularity for a single value
$t=t_0$, and the phase point $S^{t_0}x_0$ does not lie on the intersection
of more than one singularity manifolds. (In this case we say that the 
trajectory of $x_0$ has a ``simple singularity''.)

\bigskip

\subheading{Neutral Subspaces, Advance, and Sufficiency}
Consider a {\bf nonsingular} trajectory segment $S^{[a,b]}x$.
Suppose that $a$ and $b$ are {\bf not moments of collision}.

\medskip

\proclaim{Definition 2.1} The neutral space $\Cal N_0(S^{[a,b]}x)$
of the trajectory segment $S^{[a,b]}x$ at time zero ($a<0<b$)\ is
defined by the following formula:

$$
\aligned
\Cal N_0(S^{[a,b]}x)=\big \{W\in\Cal Z\colon\;\exists (\delta>0) \;
\text{s. t.} \; \forall \alpha \in (-\delta,\delta) \\
V\left(S^a\left(Q(x)+\alpha W,V(x)\right)\right)=V(S^ax)\text{ and }
V\left(S^b\left(Q(x)+\alpha W,V(x)\right)\right)=V(S^bx)\big\}.
\endaligned
$$
\endproclaim \endrem
($\Cal Z$ is the common tangent space $\Cal T_q\bold Q$ of the parallelizable
manifold $\bold Q$ at any of its points $q$, while $V(x)$ is the velocity
component of the phase point $x=\left(Q(x),\,V(x)\right)$.)

It is known (see (3) in \S 3 of [S-Ch (1987)]) that
$\Cal N_0(S^{[a,b]}x)$ is a linear subspace of $\Cal Z$ indeed, and
$V(x)\in \Cal N_0(S^{[a,b]}x)$. The neutral space $\Cal N_t(S^{[a,b]}x)$
of the segment $S^{[a,b]}x$ at time $t\in [a,b]$ is defined as follows:

$$
\Cal N_t(S^{[a,b]}x)=\Cal N_0\left(S^{[a-t,b-t]}(S^tx)\right).
$$
It is clear that the neutral space $\Cal N_t(S^{[a,b]}x)$ can be
canonically
identified with $\Cal N_0(S^{[a,b]}x)$ by the usual identification of the
tangent spaces of $\bold Q$ along the trajectory $S^{(-\infty,\infty)}x$
(see, for instance, \S 2 of [K-S-Sz(1990)]).

Our next  definition is  that of  the {\bf advance}. Consider a
non-singular orbit segment $S^{[a,b]}x$ with the symbolic collision sequence
$\Sigma=(\sigma_1, \dots, \sigma_n)$ ($n\ge 1$), meaning that $S^{[a,b]}x$
has exactly $n$ collisions with $\partial\bold Q$, and the $i$-th collision
($1\le i\le n$) takes place at the boundary of the cylinder $C_{\sigma_i}$.
For $x=(Q,V)\in\bold M$ and $W\in\Cal Z$, $\Vert W\Vert$ sufficiently small, 
denote $T_W(Q,V):=(Q+W,V)$.

\proclaim{Definition 2.2}
For any $1\le k\le n$ and $t\in[a,b]$, the advance
$$
\alpha(\sigma_k)\colon\;\Cal N_t(S^{[a,b]}x) \rightarrow  \Bbb R
$$
of the collision $\sigma_k$ is the unique linear extension of the linear 
functional $\alpha(\sigma_k)$
defined in a sufficiently small neighborhood of the origin of 
$\Cal N_t(S^{[a,b]}x)$ in the following way:
$$
\alpha(\sigma_k)(W):= t_k(x)-t_k(S^{-t}T_WS^tx).
$$
\endproclaim \endrem
Here $t_k=t_k(x)$ is the time moment of the $k$-th collision $\sigma_k$ on
the trajectory of $x$ after time $t=a$. The above formula and the notion of
the advance functional 

$$
\alpha_k=\alpha(\sigma_k):\; \Cal N_t\left(S^{[a,b]}x\right)\to\Bbb R
$$
has two important features:

\medskip

(i) If the spatial translation $(Q,V)\mapsto(Q+W,V)$ is carried out at time
$t$, then $t_k$ changes linearly in $W$, and it takes place just 
$\alpha_k(W)$ units of time earlier. (This is why it is called ``advance''.)

\medskip

(ii) If the considered reference time $t$ is somewhere between $t_{k-1}$
and $t_k$, then the neutrality of $W$ precisely means that

$$
W-\alpha_k(W)\cdot V(x)\in A_{\sigma_k},
$$
i. e. a neutral (with respect to the collision $\sigma_k$) spatial translation
$W$ with the advance $\alpha_k(W)=0$ means that the vector $W$ belongs to the
generator space $A_{\sigma_k}$ of the cylinder $C_{\sigma_k}$.

It is now time to bring up the basic notion of {\bf sufficiency} 
(or, sometimes it is also called {\bf geometric hyperbolicity}) of a
trajectory (segment). This is the utmost important necessary condition for
the proof of the fundamental theorem for semi-dispersive billiards, see
Condition (ii) of Theorem 3.6 and Definition 2.12 in [K-S-Sz(1990)].

\medskip

\proclaim{Definition 2.3}
\roster
\item
The nonsingular trajectory segment $S^{[a,b]}x$ ($a$ and $b$ are supposed not
to be moments of collision) is said to be {\bf sufficient} if and only if
the dimension of $\Cal N_t(S^{[a,b]}x)$ ($t\in [a,b]$) is minimal, i.e.
$\text{dim}\ \Cal N_t(S^{[a,b]}x)=1$.
\item
The trajectory segment $S^{[a,b]}x$ containing exactly one singularity (a so 
called ``simple singularity'', see above) is said to be {\bf sufficient} if 
and only if both branches of this trajectory segment are sufficient.
\endroster
\endproclaim \endrem

\medskip

\proclaim{Definition 2.4}
The phase point $x\in\bold M$ with at most one singularity is said
to be sufficient if and only if its whole trajectory $S^{(-\infty,\infty)}x$
is sufficient, which means, by definition, that some of its bounded
segments $S^{[a,b]}x$ are sufficient.
\endproclaim \endrem

In the case of an orbit $S^{(-\infty,\infty)}x$ with a simple
singularity, sufficiency means that both branches of
$S^{(-\infty,\infty)}x$ are sufficient.

\bigskip

\subheading{No accumulation (of collisions) in finite time} 
By the results of Vaserstein [V(1979)], Galperin [G(1981)] and
Burago-Ferleger-Kononenko [B-F-K(1998)], in a 
semi-dis\-per\-sive billiard flow there can only be finitely many 
collisions in finite time intervals, see Theorem 1 in [B-F-K(1998)]. 
Thus, the dynamics is well defined as long as the trajectory does not hit 
more than one boundary components at the same time.

\bigskip

\subheading{Slim sets} 
We are going to summarize the basic properties of codimension-two subsets $A$
of a smooth manifold $M$. Since these subsets $A$ are just those
negligible in our dynamical discussions, we shall call them {\bf slim}. 
As to a  broader exposition of the issues, see [E(1978)] or \S 2 of
[K-S-Sz(1991)].

Note that the dimension $\dim A$ of a separable metric space $A$ is one of the
three classical notions of topological dimension: the covering 
(\v Cech-Lebesgue), the small inductive (Menger-Urysohn), or the large 
inductive (Brouwer-\v Cech) dimension. As it is known from general general 
topology, all of them are the same for separable metric spaces.

\medskip

\proclaim{Definition 2.5}
A subset $A$ of $M$ is called slim if and only if $A$ can be covered by a 
countable family of codimension-two (i. e. at least two) closed sets of
$\mu$--measure zero, where $\mu$ is a smooth measure on $M$. (Cf.
Definition 2.12 of [K-S-Sz(1991)].)
\endproclaim \endrem

\medskip

\proclaim{Property 2.6} The  collection of all slim subsets of $M$ is a
$\sigma$-ideal, that is, countable unions of slim sets and arbitrary
subsets of slim sets are also slim.
\endproclaim

\medskip

\proclaim{Proposition 2.7 (Locality)}
A subset $A\subset M$ is slim if and only if for
every $x\in A$ there exists an open neighborhood $U$ of $x$ in $M$ such that 
$U\cap A$ is slim. (Cf. Lemma 2.14 of [K-S-Sz(1991)].)
\endproclaim

\medskip

\proclaim{Property 2.8} A closed subset $A\subset M$ is slim if and only
if $\mu(A)=0$ and $\dim A\le\dim M-2$.
\endproclaim

\medskip

\proclaim{Property 2.9 (Integrability)}
If $A\subset M_1\times M_2$ is a closed subset of the product of two manifolds,
and for every $x\in M_1$ the set
$$
A_x=\{ y\in M_2\colon\; (x,y)\in A\}
$$
is slim in $M_2$, then $A$ is slim in $M_1\times M_2$.
\endproclaim

\medskip

The following propositions characterize the codimension-one and 
codimension-two sets.

\proclaim{Proposition 2.10}
For any closed subset $S\subset M$ the following three conditions are 
equivalent:

\roster

\item"{(i)}" $\dim S\le\dim M-2$;

\item"{(ii)}"  $\text{int}S=\emptyset$ and for every open connected set 
$G\subset M$ the difference set $G\setminus S$ is also connected;

\item"{(iii)}" $\text{int}S=\emptyset$ and for every point $x\in M$ and for any
open neighborhood $V$ of $x$ in $M$ there exists a smaller open neighborhood
$W\subset V$ of the point $x$ such that for every pair of points 
$y,z\in W\setminus S$ there is a continuous curve $\gamma$ in the set 
$V\setminus S$ connecting the points $y$ and $z$.

\endroster

\endproclaim

\noindent
(See Theorem 1.8.13 and Problem 1.8.E of [E(1978)].)

\medskip

\proclaim{Proposition 2.11} For any subset $S\subset M$ the condition 
$\dim S\le\dim M-1$ is equivalent to $\text{int}S=\emptyset$.
(See Theorem 1.8.10 of [E(1978)].)
\endproclaim

\medskip

We recall an elementary, but important lemma (Lemma 4.15 of [K-S-Sz(1991)]).
Let $R_2$ be the set of phase points 
$x\in\bold M\setminus\partial\bold M$ such that the trajectory 
$S^{(-\infty,\infty)}x$ has more than one singularities.

\proclaim{Proposition 2.12}
The set $R_2$ is a countable union of codimension-two
smooth sub-manifolds of $M$ and, being such, it is slim.
\endproclaim

\medskip

The next lemma establishes the most important property of slim sets which
gives us the fundamental geometric tool to connect the open ergodic components
of billiard flows.

\proclaim{Proposition 2.13}
If $M$ is connected, then the complement $M\setminus A$
of a slim set $A\subset M$ necessarily contains an arc-wise connected,
$G_\delta$ set of full measure. (See Property 3 of \S 4.1 in [K-S-Sz(1989)].
The  $G_\delta$ sets are, by definition, the countable intersections
of open sets.)
\endproclaim

\medskip

\subheading{\bf The subsets $\bold M^0$ and $\bold M^\#$} Denote by
$\bold M^\#$ the set of all phase points $x\in\bold M$ for which the
trajectory of $x$ encounters infinitely many non-tangential collisions
in both time directions. The trajectories of the points 
$x\in\bold M\setminus\bold M^\#$ are lines: the motion is linear and uniform,
see the appendix of [Sz(1994)]. It is proven in lemmas A.2.1 and A.2.2
of [Sz(1994)] that the closed set $\bold M\setminus\bold M^\#$ is a finite
union of hyperplanes. It is also proven in [Sz(1994)] that, locally, the two
sides of a hyperplanar component of $\bold M\setminus\bold M^\#$ can be
connected by a positively measured beam of trajectories, hence, from the point
of view of ergodicity, in this paper it is enough to show that the connected
components of $\bold M^\#$ entirely belong to one ergodic component. This is
what we are going to do in this paper.

Denote by $\bold M^0$ the set of all phase points $x\in\bold M^\#$ the 
trajectory of which does not hit any singularity, and use the notation
$\bold M^1$ for the set of all phase points $x\in\bold M^\#$ whose orbit
contains exactly one, simple singularity. According to Proposition 2.12,
the set $\bold M^\#\setminus(\bold M^0\cup\bold M^1)$ is a countable union of
smooth, codimension-two ($\ge2$) submanifolds of $\bold M$, and, therefore,
this set may be discarded in our study of ergodicity, please see also the
properties of slim sets above. Thus, we will restrict our attention to the
phase points $x\in\bold M^0\cup\bold M^1$.

\medskip

\subheading{\bf The ``Chernov-Sinai Ansatz''} An essential precondition for
the Theorem on Local Ergodicity by B\'alint--Chernov--Sz\'asz--T\'oth is the
so called ``Chernov-Sinai Ansatz'' which we are going to formulate below.
Denote by $\Cal S\Cal R^+\subset\partial\bold M$ the set of all phase points
$x_0=(q_0,v_0)\in\partial\bold M$ corresponding to singular reflections
(a tangential or a double collision at time zero) supplied with the 
post-collision (outgoing) velocity $v_0$. It is well known that
$\Cal S\Cal R^+$ is a compact cell complex with dimension
$2d-3=\text{dim}\bold M-2$. It is also known (see Lemma 4.1 in [K-S-Sz(1990)])
that for $\nu$-almost every phase point $x_0\in\Cal S\Cal R^+$ (Here $\nu$
is the Riemannian volume of $\Cal S\Cal R^+$ induced by the restriction of
the natural Riemannian metric of $\bold M$.) the forward orbit 
$S^{(0,\infty)}x_0$ does not hit any further singularity. The Chernov-Sinai
Ansatz postulates that for $\nu$-almost every $x_0\in\Cal S\Cal R^+$
the forward orbit $S^{(0,\infty)}x_0$ is sufficient (geometrically
hyperbolic). 

\medskip

\subheading{\bf The Theorem on Local Ergodicity} The Theorem on Local 
Ergodicity by Chernov and Sinai (Theorem 5 of [S-Ch(1987)], see also
Theorem 4.4 in [B-Ch-Sz-T(2002)])
claims the following: Let $\flow$ be a semi-dispersive
billiard flow with the properties (1.1)--(1.2) and such that the smooth 
components of the boundary $\partial\bold Q$ of the configuration space 
are algebraic hypersurfaces. (The cylindric billiards with (1.1)--(1.2)
automatically fulfill this algebraicity condition.) Assume -- further --
that the Chernov-Sinai Ansatz holds true, and a phase point 
$x_0\in\bold M\setminus\partial\bold M$ is given with the properties

\medskip

(i) $S^{(-\infty,\infty)}x$ has at most one singularity,

\noindent
and

(ii) $S^{(-\infty,\infty)}x$ is sufficient. (In the case of a singular obit
$S^{(-\infty,\infty)}x$ this means that both branches of 
$S^{(-\infty,\infty)}x$ are sufficient.)

\medskip

Then some open neighborhood $U_0\subset\bold M$ of $x_0$ belongs to a single
ergodic component of the flow $\flow$. (Modulo the zero sets, of course.)

\medskip

\noindent {\bf 3. Geometric Considerations}
Consider a non-singular trajectory segment 

$$
S^{[a,b]}x_0=\left\{x_t=S^tx_0\big|\; a\le t\le b\right\}
$$
of the cylindric 
billiard flow $\flow$ with the symbolic collision sequence $\symb$, meaning
that there are time moments $a<t_1<t_2<\dots<t_n<b$ such that 
$S^tx\not\in\partial\bold M$ for $t\in[a,b]\setminus\{t_1,\dots,t_n\}$,
and $Q\left(S^{t_i}x\right)\in\partial C_{\sigma_i}$, $i=1,\dots,n$.
We assume that

\medskip

(1) $\text{dim}\left(L_{\sigma_i}\cap L_{\sigma_j}\right)\ge2$ for 
$1\le i,\, j\le n$

\noindent
(the so called ``codimension-two condition'' imposed on $\symb$), and

\medskip

(2) $\text{span}\left\{L_{\sigma_1},\dots,L_{\sigma_n}\right\}=\Bbb R^d$, 
i. e. the system of cylinders $C_{\sigma_1},\dots,C_{\sigma_n}$ is 
transitive, see also \S1.

\medskip

The first result of this section is

\medskip

\subheading{\bf Proposition 3.1} Under the above conditions (1)--(2) the 
non-singular orbit segment $S^{[a,b]}x_0$ is hyperbolic (sufficient, cf. \S2)
modulo some codimension-two (i. e. at least two) submanifolds of the phase 
space.

\medskip

\subheading{\bf Proof} The proof is based upon the following, simple lemma:

\medskip

\subheading{\bf Lemma 3.2} Let $n=2$, i. e. 
$\Sigma\left(S^{[a,b]}x_0\right)=(\sigma_1,\, \sigma_2)$. Then the advance
functionals (cf. \S2) 
$\alpha_1,\, \alpha_2:\; \Cal N_0\left(S^{[a,b]}x_0\right)\to\Bbb R$ 
corresponding to $\sigma_1$ and $\sigma_2$ are the same, unless the phase 
point $x_0$ belongs to some codimension-two ($\ge2$) submanifold.

\medskip

\subheading{\bf Proof} We may assume that the reference time $t=0$ is between
the collisions $\sigma_1$ and $\sigma_2$, i. e. 
$t_1=t(\sigma_1)<0<t_2=t(\sigma_2)$. Consider an arbitrary neutral vector
$\delta q\in\Cal N_0\left(S^{[a,b]}x_0\right)$. The neutrality of $\delta q$ 
with the advances $\alpha_i=\alpha_i(\delta q)$ ($i=1,2$) means that

$$
\delta q-\alpha_i v_0\in A_{\sigma_i}\quad (i=1,2),
\tag 3.3
$$
where $v_0=v(x_0)$, $x_0=(q_0,\, v_0)$. If $\alpha_1=\alpha_1(\delta q)$ 
happens to be different from $\alpha_2=\alpha_2(\delta q)$, then the equations
in (3.3) yield that $(\alpha_1-\alpha_2)v_0\in\text{span}\left\{A_{\sigma_1},
\, A_{\sigma_2}\right\}$, i. e. 
$v_0\in\text{span}\left\{A_{\sigma_1},\, A_{\sigma_2}\right\}$. However, 

$$
c:=\text{codim}\left(\text{span}\left\{A_{\sigma_1},\, A_{\sigma_2}\right\}
\right)=\text{dim}\left(L_{\sigma_1}\cap L_{\sigma_2}\right)\ge2
$$
(by our assumption (1)), and the event 
$v_0\in\text{span}\left\{A_{\sigma_1},\, A_{\sigma_2}\right\}$ is clearly
described by a submanifold of codimension $c$. \qed

\medskip

Finishing the proof of the proposition:

\medskip

According to the lemma, apart from a codimension-two ($\ge2$) exceptional set
$E\subset\bold M$ it is true that all advance functionals 
$\alpha_i=\alpha_{\sigma_i}:\; \Cal N\left(S^{[a,b]}x_0\right)\to\Bbb R$ 
coincide. Assume that $x_0\not\in E$ and the reference time $t=0$ is chosen 
to be right before the first collision $\sigma_1$ of $\symb$. Consider an 
arbitrary neutral vector $w=\delta q\in\Cal N_0\left(S^{[a,b]}x_0\right)$. 
By replacing $w$ with $w-\alpha v_0$ ($x_0=(q_0,v_0)$, $\alpha=\alpha_i(w)$ 
is the common value of the advances $\alpha_i(w)$) we easily achieve that 
$\alpha_i(w)=0$ for $i=1,2,\dots,n$. The relation $\alpha_1(w)=0$ means that
$w\in A_{\sigma_1}$ and $DS^{t^*_1}\left((w,0)\right)=(w,0)$, where 
$t(\sigma_1)<t^*_1<t(\sigma_2)$. Similarly, $\alpha_2(w)=0$ means that
$w\in A_{\sigma_2}$ and $DS^{t^*_2}\left((w,0)\right)=(w,0)$, where 
$t(\sigma_2)<t^*_2<t(\sigma_3)$, etc. We get that

$$
w\in\bigcap_{i=1}^n A_{\sigma_i}=\text{span}\left\{L_{\sigma_i}\big|\;
1\le i\le n\right\}^\perp=\{0\},
$$
thus $w=0$. This shows that the original neutral vector was indeed a scalar 
multiple of the velocity $v_0$, so 
$\text{dim}\Cal N_0\left(S^{[a,b]}x_0\right)=1$ whenever $x_0\not\in E$. \qed

\medskip

\subheading{\bf Remark 3.4} It is clear from the above proof that without the
assumption (1) of the proposition we obtain a codimension-one exceptional 
set $E\subset \bold M$ outside of which the statement holds true. Indeed, the
overall assumption on the geometry of our cylindric billiard system is that 
$L_i\cap L_j\ne\{0\}$ for any pair of base spaces $L_i$ and $L_j$.

\bigskip \bigskip

\heading
Some Observations Concerning Codimension-one \\
Exceptional Manifolds $J\subset\bold M$
\endheading

\bigskip \bigskip

The last remaining question of this section is this: In the original set-up 
(i. e. when only $\text{dim}\left(L_{\sigma_i}\cap L_{\sigma_j}\right)\ge1$ 
is assumed in (1) of Proposition 3.1) how an enhanced version of Proposition 
3.1 excludes the existence of a codimension-one, smooth sub-manifold 
$J\subset\bold M$ separating different ergodic components of the flow $\flow$?

\medskip

Given a codimension-one, flow-invariant, smooth sub-manifold 
$J\subset\bold M$,
consider a normal vector $n_0=(z,w)$ ($\ne 0$) of $J$ at the phase point
$y\in J$, i. e. for any tangent vector 
$(\delta q,\, \delta v)\in\Cal T_y\bold M$ the relation 
$(\delta q,\, \delta v)\in\Cal T_yJ$ is true if and only if 
$\langle\delta q,z\rangle+\langle\delta v,w\rangle=0$. Here 
$\langle\, .\, ,\, .\, \rangle$ is the Euclidean inner product of the tangent
space $\Bbb R^d$ of $\Bbb T^d$ at every point $q\in\Bbb T^d$. Let us
determine first the time-evolution $n_0\longmapsto n_t$ ($t>0$) of this normal
vector as time $t$ elapses. If there is no collision on the orbit segment
$S^{[0,t]}y$, then the relationship between
$(\delta q,\, \delta v)\in\Cal T_y\bold M$ and
$(\delta q',\, \delta v')=\left(DS^t\right)(\delta q,\, \delta v)$ is
obviously

$$
\aligned
\delta v'&=\delta v, \\
\delta q'&=\delta q+t\delta v,
\endaligned
\tag 3.5
$$
from which we obtain that

$$
\aligned
(\delta q',\, \delta v')\in\Cal T_{y'}J&\Leftrightarrow\langle\delta q'-t
\delta v',\, z\rangle+\langle\delta v',\, w\rangle=0 \\
&\Leftrightarrow
\langle\delta q',\, z\rangle+\langle\delta v',\, w-tz\rangle=0.
\endaligned
$$
This means that $n_t=(z,\, w-tz)$. It is always very useful to consider
the quadratic form $Q(n)=Q((z,w))=:\langle z,w\rangle$ associated with the
normal vector $n=(z,w)\in\Cal T_y\bold M$ of $J$ at $y$. $Q(n)$ is the
so called ``infinitesimal Lyapunov function'', see [K-B(1994)] or part
A.4 of the Appendix in [Ch(1994)]. For a detailed exposition of the 
relationship between the quadratic form $Q$, the relevant symplectic geometry
and the dynamics, please see [L-W(1995)].

\medskip

\subheading{\bf Remark} Since the normal vector $n=(z,w)$ of $J$ is only
determined up to a nonzero scalar multiplier, the value $Q(n)$ is only
determined up to a positive multiplier. However, this means that the sign
of $Q(n)$ (which is the utmost important thing for us) is uniquely
determined. This remark will gain a particular importance in the near
future.

\medskip

\noindent
From the above calculations we get that 

$$
Q(n_t)=Q(n_0)-t||z||^2\le Q(n_0).
\tag 3.6
$$

The next question is how the normal vector $n$ of $J$ gets transformed
$n^-\mapsto n^+$ through a collision (reflection) at time $t=0$? Elementary
geometric considerations show (see Lemma 2 of [Sin(1979)], or formula
(2) in \S3 of [S-Ch(1987)]) that the linearization of the flow

$$
\left(DS^t\right)\Big|_{t=0}:\; (\delta q^-,\, \delta v^-)\longmapsto
(\delta q^+,\, \delta v^+)
$$
is given by the formulae

$$
\aligned
\delta q^+&=R\delta q^-, \\
\delta v^+&=R\delta v^-+2\cos\phi RV^*KV\delta q^-,
\endaligned
\tag 3.7
$$
where the operator $R:\; \Cal T_q\bold Q\to \Cal T_q\bold Q$ is the orthogonal
reflection across the tangent hyperplane
$\Cal T_q\partial\bold Q$ of $\partial\bold Q$ at $q\in \partial\bold Q$
($y^-=(q,v^-)\in\partial\bold M$, $y^+=(q,v^+)\in\partial\bold M$), 
$V:\; (v^-)^\perp\to\Cal T_q\partial\bold Q$ is the $v^-$-parallel projection
of the orthocomplement hyperplane $(v^-)^\perp$ onto
$\Cal T_q\partial\bold Q$, $V^*:\; \Cal T_q\partial\bold Q\to (v^-)^\perp$
is the adjoint of $V$, i. e. it is the projection of $\Cal T_q\partial\bold Q$
onto $(v^-)^\perp$ being parallel to the normal vector $\nu(q)$ of
$\partial\bold Q$ at $q\in\partial\bold Q$, 
$K:\; \Cal T_q\partial\bold Q\to \Cal T_q\partial\bold Q$ is the second
fundamental form of $\partial\bold Q$ at $q$ and, finally, 
$\cos\phi=\langle\nu(q),\, v^+\rangle$ is the cosine of the angle $\phi$
subtended by $v^+$ and the normal vector $\nu(q)$. For the formula (3.7),
please also see the last displayed formula of \S1 in [S-Ch(1982)], or
(i) and (ii) of Proposition 2.3 in [K-S-Sz(1990)]. We note that it is enough
to deal with the tangent vectors 
$(\delta q^-,\, \delta v^-)\in(v^-)^\perp\times(v^-)^\perp$
($(\delta q^+,\, \delta v^+)\in(v^+)^\perp\times(v^+)^\perp$), for the
manifold $J$ under investigation is supposed to be flow-invariant, so any
vector $(\delta q,\, \delta v)=(\alpha v,\, 0)$ ($\alpha\in\Bbb R$) is
automatically inside $\Cal T_yJ$. The backward version (inverse) 

$$
\left(DS^t\right)\Big|_{t=0}:\; (\delta q^+,\, \delta v^+)\mapsto
(\delta q^-,\, \delta v^-)
$$
can be deduced easily from (3.7):

$$
\aligned
\delta q^-&=R\delta q^+, \\
\delta v^-&=R\delta v^+-2\cos\phi RV_1^*KV_1\delta q^+,
\endaligned
\tag 3.8
$$
where $V_1:\; (v^+)^\perp\to\Cal T_q\partial\bold Q$ is the $v^+$-parallel
projection of $(v^+)^\perp$ onto $\Cal T_q\partial\bold Q$. By using formula
(3.8), one easily computes the time-evolution $n^-\longmapsto n^+$
of a normal vector $n^-=(z,w)\in\Cal T_{y^-}\bold M$ of $J$ if a collision
$y^-\longmapsto y^+$ takes place at time $t=0$:

$$
\aligned
(\delta q^+,\, \delta v^+)\in\Cal T_{y^+}J\Leftrightarrow\langle R\delta q^+,
\, z\rangle+\langle R\delta v^+-2\cos\phi RV_1^*KV_1\delta q^+,\,
w\rangle &=0 \\
\Leftrightarrow\langle\delta q^+,\, Rz-2\cos\phi V_1^*KV_1Rw\rangle+
\langle\delta v^+,\, Rw\rangle &=0.
\endaligned
$$
This means that

$$
n^+=\left(Rz-2\cos\phi V_1^*KV_1Rw,\, Rw\right)
\tag 3.9
$$
if $n^-=(z,\, w)$. It follows that

$$
\aligned
Q(n^+)&=Q(n^-)-2\cos\phi\langle V_1^*KV_1Rw,\, Rw\rangle \\
&=Q(n^-)-2\cos\phi\langle KV_1Rw,\, V_1Rw\rangle\le Q(n^-).
\endaligned
\tag 3.10
$$
Here we used the fact that the second fundamental form $K$ of 
$\partial\bold Q$ at $q$ is positive semi-definite, which just means that the
billiard system is semi-dispersive.

The last simple observation on the quadratic form $Q(n)$
regards the involution $I:\; \bold M\to\bold M$, $I(q,v)=(q,-v)$ 
corresponding to the time reversal. If $n=(z,w)$ is a normal vector of $J$
at $y$, then, obviously, $I(n)=(z,-w)$ is a normal vector of $I(J)$ at
$I(y)$ and

$$
Q\left(I(n)\right)=-Q(n).
\tag 3.11
$$

By switching --- if necessary --- from the separating manifold $J$ to 
$I(J)$, and by taking a suitable remote image $S^t(J)$ ($t>>1$), in the
spirit of (3.6), (3.10)--(3.11) we can assume that

$$
Q(n)<0
\tag 3.12
$$
for every {\it unit} normal vector $n\in\Cal T_y\bold M$ of
$J$ near a phase point $y\in J$.

\medskip

\subheading{\bf Remark 3.13} There could be, however, a little difficulty in
achieving the inequality $Q(n)<0$, i. e. (3.12). Namely, it may happen that
$Q(n_t)=0$ for every $t\in\Bbb R$. According to (3.6), the equation $Q(n_t)=0$
($\forall\, t\in\Bbb R$) implies that 
$n_t=:(z_t,\, w_t)=(0,\, w_t)$ for all $t\in\Bbb R$ and, moreover, in the view
of (3.9), $w_t^+=Rw_t^-$ is the transformation law at any collision
$y_t=(q_t,\, v_t)\in\partial\bold M$. Furthermore, at every collision
$y_t=(q_t,\, v_t)\in\partial\bold M$ the projected tangent vector
$V_1Rw_t^-=V_1w_t^+$ lies in the null space of the operator $K$ 
(see also (3.9)), and this means that $w_0$ is a neutral vector for the
entire trajectory $S^{\Bbb R}y$, i. e. $w_0\in\Cal N\left(S^{\Bbb R}y\right)$.
(For the notion of neutral vectors and $\Cal N\left(S^{\Bbb R}y\right)$,
cf. \S2 above.) On the other hand, this is impossible
for the following reason: Any tangent vector $(\delta q,\delta v)$ from the
space $\Cal N\left(S^{\Bbb R}y\right)\times\Cal N\left(S^{\Bbb R}y\right)$
is automatically tangent to the separating manifold $J$, thus for any normal 
vector $n=(z,w)\in\Cal T_y\bold M$ of a separating manifold $J$ one has

$$
(z,\, w)\in\Cal N\left(S^{\Bbb R}y\right)^\perp\times\Cal N\left(
S^{\Bbb R}y\right)^\perp.
\tag 3.14
$$
(As a direct inspection shows. We always tacitly assume that the exceptional
manifold $J$ is locally defined by the equation 
$J=\left\{x\in U_0\big|\;\text{ dim}\Cal N\left(S^{[a,b]}x\right)>1\right\}$
with orbit segments $S^{[a,b]}x$ whose symbolic sequence is combinatorially
rich, i. e. it typically provides sufficient phase points.)
The membership in (3.14) is, however, impossible with a nonzero vector
$w\in\Cal N\left(S^{\Bbb R}y\right)$.

\medskip

\subheading{\bf Singularities} 

\medskip

Consider a smooth, connected piece
$\Cal S\subset\bold M$ of a singularity manifold corresponding to a
singular (tangential or double) reflection {\it in the past}. Such a
manifold $\Cal S$ is locally
flow-invariant and has one codimension, so we can
speak about its normal vectors $n$ and the uniquely determined sign of
$Q(n)$ for $0\ne n\in\Cal T_y\bold M$, $y\in\Cal S$, $n\perp\Cal S$
(depending on the foot point, of course). Consider first a phase point
$y^+\in\partial\bold M$ right after the singular reflection that is
described by $\Cal S$. It follows from the proof of Lemma 4.1 of
[K-S-Sz(1990)] and Sub-lemma 4.4 therein that at 
$y^+=(q,\, v^+)\in\partial\bold M$ any tangent vector 
$(0,\, \delta v)\in\Cal T_{y^+}\bold M$ lies actually in 
$\Cal T_{y^+}\Cal S$ and, consequently, the normal vector
$n=(z,w)\in\Cal T_{y^+}\bold M$ of $\Cal S$ at $y^+$ necessarily has the
form $n=(z,0)$, i. e. $w=0$. Thus $Q(n)=0$ for any normal vector
$n\in\Cal T_{y^+}\bold M$ of $\Cal S$. According to the monotonicity 
inequalities (3.6) and (3.10) above,

$$
Q(n)<0
\tag 3.15
$$
for any phase point $y\in\Cal S$ of a past singularity manifold $\Cal S$.

\medskip

The above observations lead to the following conclusion:

\medskip

\subheading{\bf Proposition 3.16} Assume that the separating manifold
$J\subset\bold M$ ($J$ is smooth, connected, $\text{codim}(J)=1$) is
selected in such a way that $Q(n_y)<0$ for all normal vectors
$0\ne n_y\in\Cal T_y\bold M$ of $J$ at any point $y\in J$, see above. 
Suppose further that the non-singular orbit segments $S^{[a,b]}y$
($y\in B_0$, $B_0$ is a small open ball, $0<a<b$ fixed) have the common
symbolic collision sequence $\symb$ with the relaxed properties

\medskip

(1)' $\text{dim}\left(L_{\sigma_i}\cap L_{\sigma_j}\right)\ge1$
($1\le i,\, j\le n$), and

\medskip

(2)' (the same as (2)) 
$\text{span}\left\{L_{\sigma_1},\dots,L_{\sigma_n}\right\}=\Bbb R^d$.

\medskip

We claim that for almost every phase point $y\in J\cap B_0$ (i. e. apart from
an algebraic variety $E'\subset J\cap B_0$ with 
$\text{dim}(E')<\text{dim}(J)$) the orbit segment $S^{[a,b]}y$ is hyperbolic
(sufficient).

\medskip

\subheading{\bf Proof} Consider the algebraic variety $E\subset B_0$ of
exceptional phase points characterized by the proof of Proposition 3.1,
i. e. let

$$
E=\left\{y\in B_0\big|\; S^{[a,b]}y \text{ is not hyperbolic}\right\}.
$$
The mentioned proof (in particular, the proof of Lemma 3.2) shows that the
only way to have a codimension-one, smooth component in the variety
$E$ is to have a submanifold defined by the relation

$$
v_t\in\text{span}\left\{A_{\sigma_i},\, A_{\sigma_{i+1}}\right\}
\tag 3.17
$$
for some $i\in\{1,2,\dots,n-1\}$, $t(\sigma_i)<t<t(\sigma_{i+1})$ ($v_t$ is 
the velocity at time $t$) with
$\text{dim}\text{span}\left\{A_{\sigma_i},\, A_{\sigma_{i+1}}\right\}=d-1$.
However, the normal vector $n_t=(z_t,w_t)\in\Cal T_{y_t}\bold M$ of the
manifold defined by (3.17) at the point $y_t=S^ty=(q_t,v_t)$ obviously
has the form $(z_t,w_t)=(0,w_t)$, thus $Q(n_t)=0$. By the assumption of this
proposition (and by the monotone non-increasing property of $Q(n_t)$ in $t$,
see the inequalities (3.6) and (3.10)) we get that any codimension-one, smooth
component of $E$ is transversal to $J$, thus proving the proposition. \qed

\medskip

In view of the inequality (3.15) (valid for past-singularity manifolds
$\Cal S$), the exceptional manifold $J\subset\bold M$ featuring Proposition
3.16 may be replaced by any past-singularity manifold $\Cal S$ without hurting
the proof of the proposition. Thus, we obtain

\medskip

\subheading{\bf Corollary 3.18} Let $\Cal S$ be a smooth component of a
past-singularity set, $y_0\in\Cal S$, $0<a<b$, and assume that the
non-singular orbit segment $S^{[a,b]}y_0$ has the symbolic collision sequence
$\Sigma=\symb$ fulfilling (1)'--(2)' of Proposition 3.16. Then there is an
open neighborhood $B_0$ of $y_0$ in $\bold M$ such that for almost every phase
point $y\in\Cal S\cap B_0$ (with respect to the induced hypersurface measure
of $\Cal S\cap B_0$) the symbolic collision sequence $S^{[a,b]}y$ is still
the same $\Sigma=\symb$, and $S^{[a,b]}y$ is hyperbolic. \qed

(Note that in this result, as usual, the phrase ``almost every'' may be 
replaced by saying that ``apart from a countable family of smooth, proper
sub-manifolds''.)

\bigskip

\heading
Yet Another Corollary of (3.6), (3.9), and (3.10)
\endheading

\bigskip

\subheading{\bf Corollary 3.19} Assume that the orbit segment $S^{[a,b]}y_0$
is not singular,
$\Sigma\left(S^{[a,b]}y_0\right)\allowmathbreak=\symb$ fulfills the relaxed
conditions (1)'--(2)' of Proposition 3.16 and, finally, there is a 
codimension-one, flow-invariant, smooth submanifold $E\ni y_0$ in $\bold M$ 
such that

\medskip

(i) $\Sigma\left(S^{[a,b]}y\right)=\Sigma\left(S^{[a,b]}y_0\right)$
for all $y\in E$, and

(ii) $S^{[a,b]}y$ is not hyperbolic for all $y\in E$.

\medskip

Let $0\ne n_t=(z_t,w_t)$ be a normal vector of $S^t(E)$ at the point
$y_t=(q_t,v_t)$. We claim that $Q(n_t)<0$ for all $t>b$.

\medskip

\subheading{\bf Proof} As we have seen before, the manifold $E$ is defined by
the relation 
$v_t\in\text{span}\left\{A_{\sigma_i},\, A_{\sigma_{i+1}}\right\}$ with some
$i\in\{1,2,\dots,n-1\}$, 
$\text{dim}\left(\text{span}\left\{A_{\sigma_i},\, A_{\sigma_{i+1}}\right\}
\right)=d-1$, $t(\sigma_i)<t<t(\sigma_{i+1})$. Thus $n_t=(0,\bar w)$, 
$0\ne\bar w\perp\text{span}\left\{A_{\sigma_i},\, A_{\sigma_{i+1}}\right\}$,
$t(\sigma_i)<t<t(\sigma_{i+1})$, meaning also that $Q(n_t)=0$. Assume, to the
contrary of the assertion of this corollary, that $Q(n_\tau)=0$ for some
$\tau>b$. Then, by the non-increasing property of $Q(n_\tau)$ in $\tau$,
there is a small $\epsilon>0$ such that $Q(n_\tau)=0$ for all $\tau$,
$t(\sigma_{i+1})<\tau<t(\sigma_{i+1})+\epsilon$. By (3.10) this means that
the vector $V_1R\bar w=V\bar w$ is in the null space of the operator $K$,
i. e. $V\bar w\in A_{\sigma_{i+1}}$. This means, in particular, that 
$\bar w\in\text{span}\left\{v_t,\, A_{\sigma_{i+1}}\right\}$ for
$t(\sigma_i)<t<t(\sigma_{i+1})$. On the other hand, 
$\text{span}\left\{v_t,\, A_{\sigma_{i+1}}\right\}\subset\text{span}
\left\{A_{\sigma_{i}},\, A_{\sigma_{i+1}}\right\}$, and
$\bar w\perp\text{span}\left\{A_{\sigma_{i}},\, A_{\sigma_{i+1}}\right\}$,
a contradiction. \qed

\medskip

An almost immediate consequence of the above corollary and the proof of 
Proposition 3.16 is 

\medskip

\subheading{\bf Corollary 3.20} Assume that the non-singular orbit segments
$S^{[a,b]}y$ ($y\in B_0$, $B_0$ is a small, open ball) have the common
symbolic collision sequence

$$
\left(\Sigma^{(1)},\, \Sigma^{(2)}\right)=\left(\sigma_1^{(1)},\dots,
\sigma_m^{(1)};\, \sigma_1^{(2)},\dots,\sigma_n^{(2)}\right)
$$
such that both $\Sigma^{(j)}$ are {\it combinatorially rich}, i. e.

$$
\text{span}\left\{L_{\sigma_1^{(1)}},\, L_{\sigma_2^{(1)}},\, \dots,\,
L_{\sigma_m^{(1)}}\right\}=
\text{span}\left\{L_{\sigma_1^{(2)}},\, L_{\sigma_2^{(2)}},\, \dots,\,
L_{\sigma_n^{(2)}}\right\}=\Bbb R^d.
$$
Then the exceptional set

$$
E=\left\{y\in B_0\big|\; S^{[a,b]}y \text{ is not hyperbolic}\right\}
$$
has codimension at least two.

\medskip

\subheading{\bf Proof} Let $E^{(j)}$ be a smooth, codimension-one exceptional
manifold for $\Sigma^{(j)}$, $j=1,\, 2$. (The word ``exceptional'' refers to
the fact that these manifolds consist of atypical phase points for which the
corresponding $\Sigma^{(j)}$-part of the orbit is not hyperbolic, despite
the assumed combinatorial richness of $\Sigma^{(j)}$.) Let
$t(\sigma_m^{(1)})<t<t(\sigma_1^{(2)})$. By the previous corollary, the
manifold $S^t\left(E^{(1)}\right)$ has a normal vector $n_t^{(1)}$
with $Q\left(n_t^{(1)}\right)<0$ at any point 
$y_t=S^ty$ ($y\in E^{(1)}$), while, by the same corollary again (applied in
reverse time), at any phase point $y_t=S^ty$, $y\in E^{(2)}$, the 
manifold $S^t\left(E^{(2)}\right)$ has a normal vector $n_t^{(2)}$
with $Q\left(n_t^{(2)}\right)>0$, i. e. $E^{(1)}$ and $E^{(2)}$ are 
transversal at any point of their intersection. This finishes the proof
of the corollary. \qed

\medskip

\noindent {\bf 4. Hyperbolicity Is Abundant The Inductive Proof}
Below we present the inductive proof of the Theorem of this paper. 
The induction will be performed with respect to the number of cylinders $k$.

Beside the ergodicity (and, therefore, the Bernoulli property, see [C-H(1996)])
and [O-W(1998)]) we will prove (and use as the induction hypothesis!) a few
technical properties listed below:

\medskip

(H1) The Chernov--Sinai Ansatz (see \S2 above) holds true for the cylindric
billiard flow $\flow$;

\medskip

(H2) There exists a slim subset $S\subset\bold M$ (see \S2 for the concept of
``slimness'') such that for all $x\in\bold M\setminus S$

\medskip 

\hskip 0.3truein
(i) $S^{(-\infty,\infty)}x$ has at most one singularity and

\hskip 0.3truein
(ii) $S^{(-\infty,\infty)}x$ is hyperbolic (in the singular case both branches
of $S^{(-\infty,\infty)}x$ are supposed to be hyperbolic, see \S2 above).

\medskip

Consequently, according to the Fundamental Theorem for
semi-dis\-persive billiards by Chernov and Sinai (Theorem 5 of [S-Ch(1987)],
see also Theorem 4.4 in [B-Ch-Sz-T(2002)])

\medskip

(H3) For every $x\in\bold M\setminus S$ the assertion of the Fundamental
Theorem holds true in some open neighborhood $U_0$ of $x$ in
$\bold M$, in particular, $x$ is a so called ``zig-zag point'', see
Definition 5.1 in [Sz(2000)]. Consequently, since the complementer set
$\bold M\setminus S$ is known to contain a connected set of full measure
(see \S2) and the open neighborhood $U_0$ of $x$ belongs to a single ergodic
component, we get that

\medskip

(H4) $\flow$ is ergodic, hence it is a Bernoulli flow by [C-H(1996)] and
[O-W(1998)].

\medskip

The above properties (H1)---(H2) will serve for us as the induction hypothesis.

\bigskip

\heading
1. The base of the induction: $k=1$
\endheading

\bigskip

In this case, necessarily, $L_1=\Bbb R^d$ and $A_1=\{0\}$, so the cylindric 
billiard system is actually a genuine, $d$-dimensional Sinai--billiard with a
single spherical scatterer which has been well known to enjoy the properties
(H1)---(H2) since the seminal work [S-Ch(1987)].

\bigskip

\heading
2. The induction step: $<k\longrightarrow k$, $k\ge2$.
\endheading

\bigskip

Let $k\ge2$, $\flow$ be a cylindric billiard flow fulfilling all the hypotheses
of our Theorem, and suppose that the induction hypotheses have been successfully
proven for every system (within the framework of the Theorem) with less than $k$
cylindric scatterers.

First we prove (H1) for $\flow$. The upcoming proof of the Chernov-Sinai Ansatz
is going to be a local argument by nature, for stating that ``the forward orbit
of almost every phase point $x$ on a past-singularity manifold is hyperbolic''
is a local assertion.

\medskip

Let $\Cal S_0\subset\bold M^{\#}\setminus\partial\bold M$ (For the definition 
of $\bold M^{\#}$, please see \S2.) be a small piece of a 
past-singularity manifold with the following properties:

\medskip

(1) $\Cal S_0$ is smooth (analytic) and diffeomorphic to $\Bbb R^{2d-2}$;

(2) For every phase point $x\in\Cal S_0$ the last collision on the backward
orbit $S^{(-\infty,0)}x$ is a singular collision taking place at time
$\tau(x)<0$ so that the collision at $S^{\tau(x)}x$ is a simple singularity,
see \S2. Consequently, the type of this singularity (see \S2) is the same for all
$x\in\Cal S_0$.

\medskip

We will measure the size of the subsets $A\subset\Cal S_0$ by using the 
hypersurface measure $\nu$ induced on $\Cal S_0$.

First of all, we restrict our attention to the subset

$$
A=\left\{x\in\Cal S_0\big|\; S^{(0,\infty)}x \text{ is non-singular}\right\}
\tag 4.1
$$
of $\Cal S_0$, being a $G_\delta$ subset of $\Cal S_0$ with full measure,
$\nu(\Cal S_0\setminus A)=0$, for the set $\Cal S_0\setminus A$ is a countable
union of smooth, proper submanifolds of $\Cal S_0$, see Lemma 4.1 in 
[K-S-Sz(1990)]. With each phase point $x\in A$ we associate the infinite
symbolic collision sequence 

$$
\Sigma(x)=\left(\sigma_1(x),\, \sigma_2(x),\dots\right)
$$
of the forward orbit $S^{(0,\infty)}x$, the set of cylinders

$$
C(x)=\left\{\sigma_n(x)|\quad n=1,2,\dots\right\},
\tag 4.2
$$
and the linear subspace $L(x)$ of $\Bbb R^d$ spanned by the base spaces of 
these cylinders:

$$
L(x)=\text{span}\left\{L_i\big|\; i\in C(x)\right\}.
\tag 4.3
$$

In view of Corollary 3.18, for $\nu$-almost every phase point $x$ of the
(open) subset

$$
A_0=\left\{x\in A\big|\; L(x)=\Bbb R^d\right\}
\tag 4.4
$$
of $A$ the forward orbit $S^{(0,\infty)}x$ is hyperbolic (sufficient). Thus,
in order to prove the Ansatz it is enough to show that $A_0=A$.

We argue by contradiction. Assume that

$$
A_1=\left\{x\in A\big|\; C(x)=\Cal C_0\right\},\quad
\nu(A_1)>0
\tag 4.5
$$
for some $\Cal C_0\subset\{1,2,\dots,k\}$ with

$$
L^*=\text{span}\left\{L_i\big|\; i\in\Cal C_0\right\}\ne\Bbb R^d.
$$
We will get a contradiction by applying (essentially) the invariant manifold
construction ideas borrowed from the proof of Theorem 6.1 of [Sim(1992-A)].
Indeed, we pay attention to the forward orbits $S^{(0,\infty)}x$ of the points
$x\in A_1$ governed solely by the sub-billiard dynamics defined by the cylinders
$C_i$ with $i\in\Cal C_0$ in $\Bbb T^d=\Bbb R^d/\Bbb Z^d$. This sub-billiard 
dynamics obviously does not fulfill the transitivity condition, for it is
invariant under all spatial translations by the elements of the subtorus
$\tilde A=A^*/\left(A^*\cap\Bbb Z^d\right)$, where $A^*=(L^*)^\perp$ is the 
orthogonal complement of $L^*$ in $\Bbb R^d$. (Recall that the subspaces
$L^*=\text{span}\left\{L_i\big|\; i\in\Cal C_0\right\}$ and $A^*$ are lattice
subspaces, as elementary linear algebra shows.) The orthogonal direct sum

$$
\tilde L\oplus\tilde A=L^*/\left(L^*\cap\Bbb Z^d\right)\oplus
A^*/\left(A^*\cap\Bbb Z^d\right)
$$
of the sub-tori $\tilde L$ and $\tilde A$ provides a finite covering of
$\Bbb T^d=\Bbb R^d/\Bbb Z^d$. Therefore, the sub-billiard dynamics
$\left\{S_{\Cal C_0}^t\right\}$ defined by the cylinders with index in $\Cal C_0$
is finitely covered by the direct product flow 
$\left\{S^t_*\times T^t_*\right\}$, where $\left\{S^t_*\right\}$ is the 
(transitive) cylindric billiard flow in the torus 
$\tilde L=L^*/\left(L^*\cap\Bbb Z^d\right)$ defined by the intersections of the
cylinders $C_i$ ($i\in\Cal C_0$) with the sub-torus $\tilde L$, while
$\left\{T^t_*\right\}$ is the almost periodic (uniform) motion in
$\tilde A=A^*/\left(A^*\cap\Bbb Z^d\right)$: More precisely, the phase point
$x=(q,v)\in\Bbb T^d\times\Bbb R^d$ ($x\in A_1$) is decomposed locally as
$q=q_1+q_2$, $q_1\in\tilde L$, $q_2\in\tilde A$, $v=v_1+v_2$, $v_1\in L^*$,
$v_2\in A^*=(L^*)^\perp$, $S^t\left((q,v)\right)=(q(t),v(t))$, 
$q(t)=q_1(t)+q_2(t)$, $v(t)=v_1(t)+v_2(t)$, 
$\left(q_1(t),v_1(t)\right)=S^t_*(q_1,v_1)$, $v_2(t)=v_2(0)=v_2$,
$q_2(t)=q_2+tv_2=q_2(0)+tv_2(0)$. We are about to construct the local, weakly
stable manifolds $\gamma^{ws}(x)$ for $x=(q,v)=(q_1+q_2,\, v_1+v_2)\in A_1$
in the following way:

$$
\aligned
\gamma^{ws}(x)=\bigg\{y=(q_1+\delta q_1+q_2+\delta q_2,\, v_1+\delta v_1+v_2)
\bigg| \\
\text{dist}\left(S^t_*(q_1,v_1),\, 
S^t_*(q_1+\delta q_1,v_1+\delta v_1)\right)\to 0 \\
\text{exp. fast as }t\to\infty,\text{ and }
||\delta q_1||+||\delta v_1||+||\delta q_2||<\epsilon_0\bigg\}.
\endaligned
\tag 4.6
$$
We see that $\gamma^{ws}(x)$ (if it exists as a manifold containing $x$ in its
interior) is indeed the weakly stable manifold of the phase point $x$ 
corresponding to the artificially defined dynamics $S^t_*\times T^t_*$ for
$t>0$. There are two important facts here:

\medskip

(A) The weakly stable manifolds $\gamma^{ws}(x)$ (yet to be constructed for 
typical $x\in A_1$) are concave, local orthogonal sub-manifolds (see the
``Invariant Manifolds'' part of \S2 in [K-S-Sz(1990)]) and, as such, they are
uniformly transversal to the manifold $\Cal S_0$, see Sub-lemma 4.2 in
[K-S-Sz(1990)];

\medskip

(B) The exponentially stable part

$$
\aligned
\gamma^{es}(x)=\bigg\{y=(q_1+\delta q_1+q_2,\, v_1+\delta v_1+v_2)\bigg| \\
\text{dist}\left(S^t_*(q_1,v_1),\,
S^t_*(q_1+\delta q_1,v_1+\delta v_1)\right)\to 0 \\
\text{exp. fast as }t\to\infty,\text{ and }
||\delta q_1||+||\delta v_1||<\epsilon_0\bigg\}
\endaligned
\tag 4.6/a
$$
of $\gamma^{ws}(x)$ ($x=(q_1+q_2,v_1+v_2)\in A_1$) is to be constructed by using
the statement of the Fundamental Theorem (Theorem 5 of [S-Ch(1987)])
for the $\Cal C_0$-sub-billiard system $\left\{S^t_*\right\}$. This statement
can be used, for the $\nu$-typical phase points $x=(q,v)=(q_1+q_2,v_1+v_2)$
of $A_1$ have the property that the $S_*$-part 
$\left\{S^t_*(q_1,v_1)\right\}$ of their forward orbit is hyperbolic with respect
to the sub-billiard system defined by the cylinders $C_i$, $i\in\Cal C_0$,
see Corollary 3.18.

\medskip

According to the above points (A) and (B), there exists a measurable subset
$A_2\subset A_1$ with $\nu(A_2)>0$ and a number $\delta_0>0$ such that for every
$x\in A_2$ the manifold $\gamma^{ws}(x)$ exists and its boundary is at least at
the distance $\delta_0$ from $x$ (these distances are now measured by using the
induced Riemannian metric on $\gamma^{ws}(x)$). Then, by the absolute continuity
of the foliation, see Theorem 4.1 in [K-S(1986)], the union

$$
B_2=\bigcup_{x\in A_2}\gamma^{ws}(x)\subset\bold M
$$
has a positive $\mu$-measure in the phase space $\bold M$. 

Finally, the genuine forward orbits $S^{(0,\infty)}x$ of all points $x\in A_2$
avoid a fixed open ball $B_{r_0}$ of radius $r_0>0$. (For example: We may take
any open ball $B_{r_0}$ inside the interior of any avoided cylinder
$C_j$, $j\not\in\Cal C_0$.) Therefore, the forward orbit in the direct product 
dynamics $\left(S^t_*\times T^t_*\right)(y)$ of any point $y\in B_2$
($y\in\gamma^{ws}(x)$, $x\in A_2$) avoids a slightly shrunk open ball
$B_{r_0-\delta_0}$ of reduced radius $r_0-\delta_0$. However, this is clearly
impossible, for the following reason: For $y=(q_1+q_2,\, v_1+v_2)$
($q_1\in\tilde L$, $q_2\in\tilde A$, $v_1\in L^*$, $v_2\in A^*$) the $v_2$
component is left invariant by the product flow $S^t_*\times T^t_*$, and for
almost every fixed value $v_2\in A^*$ (namely, for those vectors $v_2$
for which the orbit $tv_2/\left(A^*\cap\Bbb Z^d\right)$ ($t\in\Bbb R$) is dense
in the torus $\tilde A=A^*/\left(A^*\cap\Bbb Z^d\right)$) the product flow 
$S^t_*\times T^t_*$ is ergodic on the corresponding level set, since it is the
product of a mixing and an ergodic flow. (The flow $S^t_*$ is mixing by the
induction hypothesis (H1)--(H4).) The obtained contradiction finishes
the indirect proof of the Chernov-Sinai Ansatz, that is, (H1). \qed

\bigskip

\subheading{\bf Corollary 4.7} The set

$$
NH(\Cal S_0)=\left\{x\in\Cal S_0\big|\; S^{(0,\infty)}x 
\text{ is not hyperbolic}\right\}
$$
is a slim set. (In the case of a singular forward orbit non-hyperbolicity of
$S^{(0,\infty)}x$ is meant that at least one branch of $S^{(0,\infty)}x$
is not hyperbolic, see \S2.)

\medskip

\subheading{\bf Proof} Since the complement set 
$\Cal S_0\setminus A=\Cal S_0\setminus A_0$ is a countable union of smooth,
proper sub-manifolds of $\Cal S_0$, the set $\Cal S_0\setminus A_0$ is slim.
Therefore, it is enough to prove that the intersection
$NH(\Cal S_0)\cap A_0$ is slim. However, according to Corollary 3.18,
the forward orbit $S^{(0,\infty)}x$ of every $x\in A_0$ is hyperbolic, unless
$x$ belongs to a countable union of smooth, proper sub-manifolds of $\Cal S_0$.
Thus $NH(\Cal S_0)\cap A_0$ is slim. \qed

\medskip

In view of Lemma 4.1 of [K-S-Sz(1990)], the set $R_2$ of phase points with
more than one singularity on their orbit is slim, see also \S2. Therefore, 
the final step in proving the remaining unproven induction hypothesis (i. e. 
(H2)) for our considered model $\flow$ with $k$ ($\ge2$) cylindric 
scatterers is to show that the set

$$
D=\left\{x\in\bold M^0\setminus\partial\bold M
\big|\; S^{(-\infty,\infty)}x \text{ is not hyperbolic}\right\}
\tag 4.8
$$
is slim, i. e. it can be covered by a countable collection of closed subsets
$F\subset\bold M$ with $\mu(F)=0$ and $\text{dim}F\le\text{dim}\bold M-2$.
By the locality of slimness, see \S2 above, it is enough to prove that for 
every element $x\in D$ the point $x$ has an open neighborhood $U$ (in $\bold M$)
such that the set $U\cap D$ is slim. We want to classify the phase points
$x\in D$.

Consider, therefore, an arbitrary phase point $x=(q,v)\in D$. Denote the doubly
infinite, symbolic collision sequence of $S^{(-\infty,\infty)}x$ by
$\Sigma=\left(\dots,\sigma_{-2},\sigma_{-1},\,\sigma_{1},\sigma_{2},\dots\right)$
so that $\sigma_1$ is the first collision in positive time. (The index $0$ is not
used.) We distinguish between two cases:

\medskip

\subheading{Case I}
$L^*=\text{span}\left\{L_{\sigma_i}|\; 
i\in\Bbb Z\setminus\{0\}\right\}\ne\Bbb R^d$. 

In this case, as we have seen before, the dynamics of $S^{(-\infty,\infty)}x$
is finitely covered by the direct product flow 
$\left\{S^t_*\times T^t_*\right\}$, where $\left\{S^t_*\right\}$ is the 
cylindric billiard flow in the sub-torus 
$\tilde L=L^*/\left(L^*\cap\Bbb Z^d\right)$ with
the scatterers $C_{\sigma_i}\cap\tilde L$, while $\left\{T^t_*\right\}$
is the almost periodic (uniform) motion in the orthocomplement torus
$\tilde A=A^*/\left(A^*\cap\Bbb Z^d\right)$, $A^*=(L^*)^\perp$. 
Now the point is that for the cylindric billiard flow 
$(\tilde L,\left\{S^t_*\right\},\mu_{\tilde L})$ both of the induction 
hypotheses (H1)--(H2) and, consequently, Theorem 5.2 of [Sz(2000)] apply. 
For the phase point $x\in D$ the direct product flow 
$\left(S^t_*\times T^t_*\right)(x)$ avoids an open ball, namely any open ball 
in the interior of any avoided cylinder $C_j$ with

$$
j\not\in\left\{\sigma_i|\; i\in\Bbb Z\setminus\{0\}\right\}.
$$
Consequently, for each component $(q_2,v_2)\in\tilde A\times A^*$ of the 
canonical decomposition of $x=(q,v)=(q_1+q_2,v_1+v_2)$, $q_1\in\tilde L$,
$v_1\in L^*$, $q_2\in\tilde A$, $v_2\in A^*$ it is true that the
$\tilde L$-orbit $S^t_*(q_1,v_1)$ of $(q_1,v_1)$ avoids an open set
$\emptyset\ne B\subset\tilde L$ on a doubly unbounded set $H$ of time moments,
$\inf H=-\infty$, $\sup H=+\infty$. Therefore, in view of Theorem 5.2 of
[Sz(2000)], the $(q_1,v_1)$-part of the phase point $x=(q_1+q_2,v_1+v_2)$
belongs to a slim subset $S_1$ of the phase space $\tilde L\times L^*$.
According to the integrability property of closed slim sets (cf. Property 4
in \S4.1 of [K-S-Sz(1989)]), even the closure $\bar D_1$ of the set

$$
D_1=\left\{x\in D\big|\; \text{span}\left\{L_{\sigma_i(x)}|\; i\in\Bbb Z
\setminus\{0\}\right\}\ne\Bbb R^d\right\}
\tag 4.9
$$
(covered by Case I) is a slim subset of the phase space $\bold M$. We note that
the set $\bar D_1$ is contained in the closed zero-set

$$
\aligned
K=\big\{x\in\bold M^\#\big|\; x
\text{ has a trajectory branch with a symbolic sequence } \\
(\dots,\, \sigma_{-1},\, \sigma_1,\, \dots)
\text{ such that }\text{span}\left\{L_{\sigma_i}|\; i\in\Bbb Z\setminus\{0\}
\right\}\ne\Bbb R^d\big\},
\endaligned
$$
and the argument with ``integrating up'' the closed slim sets (by using
Property 4 in \S4.1 of [K-S-Sz(1989)]) is applied to the closed set $K$.

\bigskip

\subheading{Case II} $L^*=\text{span}\left\{L_{\sigma_i}|\; 
i\in\Bbb Z\setminus\{0\}\right\}=\Bbb R^d$. 

Select a vector $0\ne w\in\Cal N\left(S^{(-\infty,\infty)}x\right)$, 
$w\perp v$, from the neutral space

$$
\Cal N\left(S^{(-\infty,\infty)}x\right)=\Cal N(x)
$$ 
of the considered phase 
point $x\in D$. For $i\in\Bbb Z\setminus\{0\}$ denote by $\alpha_i=\alpha_i(w)$
the ``advance'' of the collision $\sigma_i$ corresponding to the neutral vector
$w$, see \S2. Since $w$ is not parallel to $v$, at least two advances with
neighboring indices are unequal; we may assume that $\alpha_{-1}\ne\alpha_1$.

It follows from the proof of Lemma 3.2 that the event $\alpha_{-1}\ne\alpha_1$
can only occur if

$$
v=v_0\in\text{span}\left\{A_{\sigma_{-1}},\, A_{\sigma_{1}}\right\}.
\tag 4.10
$$
If the event $\alpha_k\ne\alpha_{k+1}$ ($k\ne-1,\, 0$) took place for another
pair of neighboring advances as well, then, again by the proof of Lemma 3.2,
we would have

$$
v_t\in\text{span}\left\{A_{\sigma_{k}},\, A_{\sigma_{k+1}}\right\}
\quad (t_k<t<t_{k+1}).
\tag 4.11
$$
If at least one of the two subspaces on the right-hand-sides of (4.10) and (4.11)
is of codimension higher than one, then the corresponding event alone
ensures that the studied phase point $x\in D$ belongs to some codimension-two
(i. e. at least two), smooth submanifold of the phase space, and such phase 
points obviously constitute a slim set, therefore they may be discarded.

Thus, we may assume that

$$
\aligned
\text{dim}\left(\text{span}\left\{A_{\sigma_{-1}},\, A_{\sigma_{1}}
\right\}\right) \\
=\text{dim}\left(\text{span}\left\{A_{\sigma_{k}},\, A_{\sigma_{k+1}}
\right\}\right)=d-1.
\endaligned
\tag 4.12
$$
Denote by $n_\tau$ a (unit) normal vector of the $S^\tau$-image of the 
manifold 

$$
v_0\in\text{span}\left\{A_{\sigma_{-1}},\, A_{\sigma_{1}}\right\}
$$
at the phase point $S^\tau x=x_\tau$, and by $\tilde n_\tau$ a (unit) 
normal vector of the $S^\tau$-image of the manifold 
$v_t\in\text{span}\left\{A_{\sigma_{k}},\, A_{\sigma_{k+1}}\right\}$ 
($t_k<t<t_{k+1}$) at the phase point $S^\tau x=x_\tau$. It follows from the
proof of Corollary 3.19 that $Q(n_\tau)=0$ for $t_{-1}<\tau<t_1$,
$Q(n_\tau)<0$ for $\tau>t_1$, $Q(n_\tau)>0$
for $\tau<t_{-1}$, $Q(\tilde n_\tau)=0$ for $t_k<\tau<t_{k+1}$,
$Q(\tilde n_\tau)<0$ for $\tau>t_{k+1}$, and
$Q(\tilde n_\tau)>0$ for $\tau<t_{k}$. Therefore, the two codimension-one
sub-manifolds defined by (4.10) and (4.11) are transversal, so the simultaneous
validity of (4.10)--(4.11) again results in an event for $x\in D$ showing that
$x$ belongs to a slim subset of $\bold M$.

\medskip

Thus we may assume that

$$
\alpha_k=\alpha_{-1}\ne\alpha_1=\alpha_l
$$
for all $k\le -1$, $l\ge1$. By adding a suitable, scalar multiple of the 
velocity $v=v_0$ to the neutral vector $w$, we can achieve that

$$
\alpha_k=\alpha_{-1}\ne0=\alpha_l
\tag 4.13
$$
for all $k\le -1$, $l\ge1$. The equalities $\alpha_l=0$ ($l\ge1$) mean that

$$
w\in\bigcap_{l>0} A_{\sigma_l}=\left(\text{span}\left\{L_{\sigma_l}\big|\;
l>0\right\}\right)^\perp,
\tag 4.14
$$
see also the closing part of the proof of Proposition 3.1. An analogous
argument shows that

$$
w-\alpha_{-1}v
\in\bigcap_{k<0} A_{\sigma_k}=\left(\text{span}\left\{L_{\sigma_k}\big|\;
k<0\right\}\right)^\perp.
\tag 4.15
$$
The equations (4.14)--(4.15) and $\alpha_{-1}\ne0$ imply that

$$
\aligned
v\in\text{span}\left\{\bigcap_{k<0}A_{\sigma_k},\,
\bigcap_{l>0}A_{\sigma_l}\right\} \\
=\bigcap_{k<0}A_{\sigma_k}+\bigcap_{l>0}A_{\sigma_l}:=H.
\endaligned
\tag 4.16
$$
Recall that $\bigcap_{n\ne0}A_{\sigma_n}=\{0\}$ in the actual Case II, and
$H\ne\Bbb R^d$, since 
$\text{span}\allowmathbreak\left\{A_{\sigma_k},\, A_{\sigma_l}\right\}\ne
\Bbb R^d$ for $k<0<l$.
We can assume that the linear direct sum on the right-hand-side of (4.16)
is a subspace with codimension one, otherwise (just as many times in the past)
the containment in (4.16) would be a codimension-two condition on the initial
velocity $v=v_0$, and all such phase points $x=(q,v)$ can be discarded. 
Following the tradition, denote by $J$ the codimension-one sub-manifold of
$\bold M$ defined by (4.16). The proof of the Ansatz in the case of (4.5)
can now be repeated almost word-by-word. Indeed, the phase points 

$$
\aligned
x=(q,\, v)\in\bar D=\bar D(\Cal A,\Cal B)=\big\{x\in D\big|\; v\in H, \\
\left\{\sigma_k(x)|\; k<0\right\}=\Cal A,\;
\left\{\sigma_l(x)|\; l>0\right\}=\Cal B\big\}
\endaligned
\tag 4.17
$$
(with given $\Cal A,\, \Cal B\subset\{1,2,\dots,k\}$ such that
$\text{span}\left\{L_j|\; j\in\Cal A\cup\Cal B\right\}=\Bbb R^d$)
of the considered type again decompose as
$(q,v)=(q_1+q_2,\, v_1+v_2)$,
$v_1\in L^*=\text{span}\allowmathbreak\left\{L_{\sigma_l}|\; l>0\right\}$, 
$v_2\in A^*=\left(L^*\right)^\perp=\bigcap_{l>0}A_{\sigma_l}$,
$q_1\in\tilde L=L^*/\left(L^*\cap\Bbb Z^d\right)$, 
$q_2\in\tilde A=A^*/\left(A^*\cap\Bbb Z^d\right)$, and the forward orbit
$S^{(0,\infty)}x$ of our considered phase point $x\in D\cap J$ (fulfilling all
of the mentioned assumptions) is essentially (up-to a finite covering) is 
governed by the product flow $\left(q_1(t),v_1(t)\right)=S^t_*(q_1,v_1)$,
$\left(q_2(t),v_2(t)\right)=T^t_*(q_2,v_2)=(q_2+tv_2,v_2)$, where (as said 
before) $S^t_*$ is the sub-billiard flow in $\tilde L$ defined by the 
intersections of the cylinders $\left\{C_{\sigma_l}|\; l>0\right\}$
with the torus $\tilde L$.

\medskip

\subheading{\bf Lemma 4.18} The exponentially stable component $\gamma^{es}(x)$
of $\gamma^{ws}(x)$ ($x\in D\cap J$) defined by (4.6/a) is transversal to the
codimension-one manifold $J$ described by the membership in (4.16).

\medskip

\subheading{\bf Proof} Argue by contradiction. Assume that
$\Cal T_x\gamma^{es}(x)\subset\Cal T_xJ$. The tangent space $\Cal T_xJ$ is
obviously given by the simple formula

$$
\Cal T_xJ=\left\{(\delta q,\delta v)\in\Cal T_x\bold M\big|\; \delta v\in H
\right\}.
\tag 4.19
$$
The second fundamental form $B\left(\gamma^{es}(x)\right)$ of $\gamma^{es}(x)$
at the phase point $x=(q_1+q_2,v_1+v_2)$ is known to be negative definite, so
its range is the entire orthocomplement $(v_1)^\perp$ of $v_1$ in the space

$$
L^*=\text{span}\left\{L_{\sigma_l}|\; l>0\right\}=
\left(\bigcap_{l>0}A_{\sigma_l}\right)^\perp.
$$
On the other hand, since 

$$
v=v_1+v_2\in H=\bigcap_{k<0}A_{\sigma_k}+\bigcap_{l>0}A_{\sigma_l}
$$
and $v_2\in A^*=\bigcap_{l>0}A_{\sigma_l}\subset H$, from the assumed relation
$\Cal T_x\gamma^{es}(x)\subset\Cal T_xJ$ and from $v_1\in H$ we get that 
$L^*\subset H$. Since $A^*=\left(L^*\right)^\perp\subset H$, this means that
$H=\Bbb R^d$, contradicting $\text{dim}H=d-1$. This finishes the proof of
the lemma. \qed

\medskip

Finally, the slimness of the set $D$ in (4.8) will be proven in Case II as soon
as we show that $\nu_J(\bar D)=0$, where $\bar D=\bar D(\Cal A,\Cal B)$
is defined in (4.17). This is, however, obtained the same way as the relation 
$\nu(A_2)=0$ at the end of the proof of the Ansatz. Indeed, in the case 
$\nu_J(\bar D)>0$ the union 

$$
\tilde D:=\bigcup_{x\in\bar D}\gamma^{es}(x)
$$
would have a positive $\mu$-measure in $\bold M$ (by the transversality proved
above and by the absolute continuity of the 
$\gamma^{es}(\,.\,)$ foliation, see Theorem 4.1 in [K-S(1986)]),
but this is impossible, for all forward orbits $S^{(0,\infty)}y$ of the points
$y\in\tilde D$ would avoid a common open ball that can be obtained by slightly
shrinking any open ball inside the interior of any avoided cylinder $C_j$ with
$j\not\in\Cal B$, see also the closing part of the proof of $\nu(A_2)=0$
above.

This finishes the proof of the fact that the set $D$ in (4.8) is indeed slim.
From this, from the proved Chernov-Sinai Ansatz, and from the quoted slimness
of the set $R_2$ of phase points with more than one singularities on their
orbit we obtain the validity of the induction hypotheses (H1)---(H2)
(and therefore (H3)---(H4), as well) for the considered cylindric billiard
flow $\flow$ with $k$ cylinders. This finishes the inductive proof of the
Theorem. \qed

\medskip

\subheading{Acknowledgement} The author expresses his sincere gratitude to
the reviewers of the paper for their careful work, especially for noticing
a few annoying mistakes in the manuscript.

E-mail address : simanyi\@math.uab.edu

\enddocument